\documentclass{amsart}%
\usepackage{amscd}
\usepackage{graphicx}
\usepackage{amsmath}
\usepackage{amsfonts}
\usepackage{amssymb}%
\newtheorem{theorem}{Theorem}
\theoremstyle{plain}

\numberwithin{equation}{section}

\begin{document}
\title{On Stochastic Generators of Positive Definite Exponents.}
\author{V. P. Belavkin.}
\address{Mathematics Department, University of Nottingham, \\
NG7 2RD, UK.}
\email{vpb@@maths.nott.ac.uk}
\thanks{This work was supported by Royal Society grant for the UK-Japan collaboration.}
\thanks{Published in: \textit{Transactions of a Japanese-German Symposium Infinite
Dimensional Harmonic Analysis} 84--92, Eds: H. Heyer et al, Kyoto, 2000.}
\date{April 2000}
\subjclass{Stochastic Analysis and Probability Theory.}
\maketitle

\begin{abstract}
A characterisation of quantum stochastic positive definite (PD) exponent is
given in terms of the conditional positive definiteness (CPD) of their
form-generator. The pseudo-Hilbert dilation of the stochastic form-generator
and the pre-Hilbert dilation of the corresponding dissipator is found. The
structure of quasi-Poisson stochastic generators giving rise to a quantum
stochastic birth processes is studied.

\end{abstract}

\section{Introduction}

Quantum probability theory provides examples of positive-definite (PD)
infinitely-divisible functions on non-Abelian groups which serve as
characteristic functions of quantum chaotic states, generalizing the
characteristic functions of classical stochastic processes with independent
increments. The simplest examples are given by quantum point processes
\cite{1} which are characterized by analytical functions on the unit ball
$B=\left\{  y\in\mathcal{B}:||y||\leq1\right\}  $ of a non-commutative group
C*-algebra. Such processes generate Markov quantum dynamics by one-parameter
families $\phi=\left(  \phi_{t}\right)  _{t>0}$ of nonlinear completely
positive maps $\phi_{t}:B\rightarrow\mathcal{A}$ on the unit ball of a
C*-algebra $\mathcal{B}$, into an operator algebra $\mathcal{A}$ of a Hilbert
space $\mathcal{H}$. As in the linear case, an analytical map $\phi_{t}$ is
completely positive iff it is positive definite (PD),
\begin{equation}
\sum_{x,z\in B}\left\langle \eta^{x}|\phi\left(  x^{\star}z\right)  \eta
^{z}\right\rangle :=\sum_{i.k}\left\langle \eta_{i}|\phi\left(  y_{i}^{\star
}y_{k}\right)  \eta_{k}\right\rangle \geq0,\quad\forall\eta_{j}\in
\mathcal{H},y_{j}\in B, \label{0.1}%
\end{equation}
where $\eta^{y}=\eta_{j}\neq0$ only for $y=y_{j},j=1,2,...$ . The simplest
quantum point dynamics of this kind is given by the quantum Markov birth
process which is described by the one-parameter semigroup
\[
\phi_{s}\left(  y\right)  \phi_{r}\left(  y\right)  =\phi_{s+r}\left(
y\right)  ,\quad\phi_{0}\left(  y\right)  =1,\quad y\in B
\]
of infinitely divisible bounded PD functions $\phi_{t}:$ $B\rightarrow
\mathbb{C}$ with the normalization property $\phi_{t}\left(  1\right)  =1$,
where $1\in B$ is (approximative) identity of $\mathcal{B}$ . The continuity
of the semigroup $\phi$ suggests the exponential form $\phi_{t}\left(
y\right)  =\exp\left[  t\lambda\left(  y\right)  \right]  $ of the functions
$\phi_{t}$ . The corresponding analytic generator
\[
\lambda\left(  y\right)  =\frac{1}{t}\ln\phi_{t}\left(  y\right)
:=\lim_{t\searrow0}\frac{1}{t}\left(  \phi_{t}\left(  y\right)  -1\right)
\]
of such semigroup is conditionally completely definite (CPD), and this is
equivalent to the PD property (\ref{0.1}) for $\phi=\lambda$ under the
condition $\sum_{j}\eta^{j}=0$ and $\lambda\left(  1\right)  =0$ . The CPD
functions have been studied in \cite{2} and the corresponding dilations
$\phi_{t}\left(  y\right)  =\left\langle \pi_{t}\left(  y\right)
\right\rangle $ to the multiplicative stochastic exponents $\pi_{t}\left(
y\right)  =:\exp\Lambda\left(  t,y\right)  :$ of a quantum process
$\Lambda\left(  t,y\right)  $ with independent increments and the vacuum mean
$\left\langle \Lambda\left(  t,y\right)  \right\rangle =t\lambda\left(
y\right)  $ in Fock space were obtained in \cite{3,4}. The unital $\star$
-multiplicative property
\[
\pi_{t}\left(  x^{\star}z\right)  =\pi_{t}\left(  x\right)  ^{\dagger}\pi
_{t}\left(  z\right)  ,\quad\pi_{t}\left(  1\right)  =I,
\]
obviously implies the PD (\ref{0.1}) of $\phi=\pi_{t}$ , and the stationarity
of the increments $\Lambda^{s}\left(  t\right)  =\Lambda\left(  t+s\right)
-\Lambda\left(  s\right)  $ implies the cocycle exponential property
\[
\pi_{s}\left(  y\right)  \pi_{r}^{s}\left(  y\right)  =\pi_{r+s}\left(
y\right)  ,\quad\forall r,s>0,
\]
with respect to the natural time-shift $\pi\mapsto\pi^{s}$ in the Fock space
of the representation $\pi$ . The dilation of the CPD generators $\lambda$
over the suggests their general form $\lambda\left(  y\right)  =\varphi\left(
y\right)  -\kappa$ , where $\varphi$ is a PD function on $B$ with
$\varphi\left(  0\right)  =0$ and $\kappa=\varphi\left(  1\right)  $ .

Here we shall extend this dilation theorem to the stochastic PD families
$\phi$ satisfying the cocycle exponential property
\[
\phi_{s}\left(  y\right)  \phi_{r}^{s}\left(  y\right)  =\phi_{r+s}\left(
y\right)  ,\quad\forall r,s>0,
\]
but not yet the unital multiplicative property. In particular, we shall obtain
the structure of the stochastic form-generator for a family $\phi$ of PD
functions $\phi_{t}\left(  \omega\right)  :B\rightarrow\mathbb{C}$ , given as
the adapted stochastic process $\phi_{t}\left(  \omega,y\right)  $ for each
$y\in B$ with respect to a classical process $\omega=\left\{  \omega\left(
t\right)  \right\}  $ with independent increments, and having the cocycle
exponential property with respect to the time-shift $\phi_{t}^{s}\left(
\omega\right)  =\phi_{t}\left(  \omega^{s}\right)  $, $\omega^{s}=\left\{
\omega\left(  t+s\right)  \right\}  $ . Such stochastic functions can be
unbounded, but they are usually normalized, $\phi_{t}\left(  \omega,1\right)
=m_{t}\left(  \omega\right)  $, to a positive-valued process $m_{t}\geq0$ ,
having the martingale property
\[
m_{t}\left(  \omega\right)  =\epsilon_{t}\left[  m_{s}\right]  \left(
\omega\right)  ,\quad\forall s>t,\quad m_{0}\left(  \omega\right)  =1,
\]
where $\epsilon_{t}$ is the conditional expectation with respect to the
history of the process $\omega$ up to time $t$ . As follows from our dilation
theorem, for example the stochastic exponent
\[
\phi_{t}\left(  y\right)  =\left(  1+\alpha\left(  y\right)  \right)
^{p\left(  t\right)  }\exp\left[  t\lambda\left(  y\right)  \right]
\]
with respect to the standard Poisson process $p\left(  t,\omega\right)  $ is
PD and normalized in the mean iff $1+\alpha$ and $\kappa+\lambda$ are PD for a
$\kappa\geq0$ , and $\alpha\left(  1\right)  +\lambda\left(  1\right)  =0.$

\section{The Generators of Quantum Stochastic PD Exponents.}

Let us consider a (noncommutative) It\^{o} $\flat$ -algebra $\mathfrak{a} $
\cite{4,0}, i.e. an associative $\star$ -algebra, identified with the algebra
of quadruples $\boldsymbol{a}=\left(  a_{\nu}^{\mu}\right)  _{\nu=+,\bullet
}^{\mu=-,\bullet} $ ,
\[
a_{\bullet}^{\bullet}=i\left(  a\right)  ,\quad a_{+}^{\bullet}=k\left(
a\right)  ,\quad a_{\bullet}^{-}=k^{*}\left(  a\right)  ,\quad a_{+}%
^{-}=l\left(  a\right)  ,
\]
under the product $\boldsymbol{b}\boldsymbol{a}=\left(  b_{\bullet}^{\mu
}a_{\nu}^{\bullet}\right)  $ and the involution $a\mapsto b=a^{\star}%
\in\mathfrak{a} $ , $b^{\star}=a $ , represented by the quadruples
$\boldsymbol{b}=\boldsymbol{a}^{\flat} $ with $b_{-\nu}^{\mu}=a_{-\mu}%
^{\nu\dagger} $ , where $-\pm=\mp$ , $-\bullet=\bullet$ . Here $i\left(
b\right)  k\left(  a\right)  =k\left(  ba\right)  $ is the GNS $\star$
-representation $i\left(  a^{\star}\right)  =i\left(  a\right)  ^{\dagger} $
associated with a linear positive $\star$ -functional $l:\mathfrak{a}%
\mapsto\mathbb{C} $ , $l\left(  a^{\star}\right)  =l\left(  a\right)  ^{*} $ ,
and $k^{*}\left(  a^{\star}\right)  =k\left(  a\right)  ^{\dagger} $ is the
linear functional on the pre-Hilbert space $\mathcal{K} $ of the Kolmogorov
decomposition $l\left(  a^{\star}a\right)  =k\left(  a\right)  ^{\dagger
}k\left(  a\right)  $ of the functional $l $ , separating $\mathfrak{a} $ in
the sense $a=0\Leftrightarrow i\left(  a\right)  =k\left(  a\right)  =l\left(
a\right)  =0 $ .

Let $B $ denote a (noncommutative) semigroup with identity $1\in B $ and
involution $y\mapsto y^{\star}\in B $ , $\left(  x^{\star}z\right)  ^{\star
}=z^{\star}x,\forall x,y,z\in B $ , say, a (noncommutative) group with
$y^{\star}=y^{-1} $ , or the unital semigroup $B=1\oplus\mathfrak{b} $ of a
$\star$ -algebra $\mathfrak{b} $ with $\left(  1\oplus a\right)  ^{\star
}\left(  1\oplus c\right)  =1\oplus a\star c $ , where $a\star c=c+a^{\star
}c+a^{\star} $ for $a,c\in\mathfrak{b} $ . The stochastically differentiable
operator-valued exponent $\phi_{t}\left(  y\right)  $ over $B $ with respect
to a quantum stationary process, with independent increments $\Lambda
^{s}\left(  t\right)  =\Lambda\left(  t+s\right)  -\Lambda\left(  s\right)  $
generated by a separable It\^{o} algebra $\mathfrak{a} $ is described by the
quantum stochastic equation
\begin{equation}
\mathrm{d}\phi_{t}\left(  y\right)  =\phi_{t}\left(  y\right)
\boldsymbol{\alpha}\left(  y\right)  \mathrm{d}\boldsymbol{A}\left(  t\right)
:=\phi_{t}\left(  y\right)  \sum_{\mu,\nu}\alpha_{\nu}^{\mu}\left(  y\right)
\mathrm{d}A_{\mu}^{\nu},\qquad\text{ }y\in B\qquad\label{1.1}%
\end{equation}
with the initial condition $\phi_{0}\left(  y\right)  =I $ , for all $y\in B $
. Here $\boldsymbol{\alpha}\left(  y\right)  \in\mathfrak{a} $ is given by the
quadruple $\alpha_{\bullet}^{\bullet}=\left[  \alpha_{n}^{m}\right]  $ ,
$\alpha_{-}^{\bullet}=\left[  \alpha_{+}^{m}\right]  $ , $\alpha_{\bullet}%
^{-}=\left[  \alpha_{n}^{-}\right]  $ , $\alpha_{+}^{-} $ of complex functions
$\alpha_{\nu}^{\mu}:B\rightarrow\mathbb{C} $ , $\mu\in\left\{
-,1,2,...\right\}  ,\quad$ $\nu\in\left\{  +,1,2,...\right\}  $ and
$\boldsymbol{A}=\left(  A_{\mu}^{\nu}\right)  _{\mu=-,\bullet}^{\nu=+,\bullet}
$ is the quadruple of the canonical integrators given by the standard time
$A_{-}^{+}\left(  t\right)  =tI $ , annihilation $A_{-}^{n}\left(  t\right)  $
, creation $A_{m}^{+}\left(  t\right)  $ and exchange $A_{m}^{n}\left(
t\right)  $ operators in Fock space over $L^{2}\left(  \mathbb{R}_{+}%
\times\mathbb{N}\right)  $ with $m,n\in\mathbb{N=}\left\{  1,2,...\right\}  $
. The infinitesimal increments $\mathrm{d}A_{\mu}^{\nu}=A_{\mu}^{t\nu}\left(
\mathrm{d}t\right)  $ are formally defined by the Hudson-Parthasarathy
multiplication table \cite{5} and the $\flat$ -property \cite{4},
\begin{equation}
\mathrm{d}A_{\mu}^{\beta}\mathrm{d}A_{\gamma}^{\nu}=\delta_{\gamma}^{\beta
}\mathrm{d}A_{\mu}^{\nu},\qquad\text{ }\boldsymbol{A}^{\flat}=\boldsymbol{A}%
,\qquad\label{1.2}%
\end{equation}
where $\delta_{\gamma}^{\beta}$ is the usual Kronecker delta restricted to the
indices $\beta\in\left\{  -,1,2,...\right\}  ,\quad\gamma\in\left\{
+,1,2,...\right\}  $ and $A_{-\nu}^{\flat\mu}=A_{-\mu}^{\nu\dagger} $ with
respect to the reflection of the indices $\left(  -,+\right)  $ only. The
structural functions $\alpha_{\nu}^{\mu}$ for the $* $ -cocycles $\phi_{t}%
^{*}=\phi_{t} $ , where $\phi_{t}^{*}\left(  y\right)  =\phi_{t}\left(
y^{\star}\right)  ^{\dagger} $ should obviously satisfy the $\flat$ -property
$\boldsymbol{\alpha}^{\flat}=\boldsymbol{\alpha}$ , where $\alpha_{-\mu
}^{\flat\nu}=\alpha_{-\nu}^{\mu*} $ , $\alpha_{\nu}^{\mu*}\left(  y\right)
=\alpha_{\nu}^{\mu}\left(  y^{\star}\right)  ^{\dagger} $ even in the case of
nonlinear $\alpha_{\nu}^{\mu}$ . The summation in (\ref{1.1}) is defined as a
quantum stochastic differential \cite{4} if $\sum_{n=1}^{\infty}\alpha_{n}%
^{-}\left(  y^{\star}\right)  \alpha_{+}^{n}\left(  y\right)  <\infty$ and the
matrix $\left[  \alpha_{n}^{m}\left(  y\right)  \right]  ,m,n\in\mathbb{N} $
represents a bounded operator in the Hilbert space $\mathbb{\ell}_{\mathbb{N}%
}^{2}=\left\{  \zeta^{\bullet}:\mathbb{N\rightarrow C}|\sum_{n=1}^{\infty
}|\zeta^{n}|^{2}<\infty\right\}  $ for each $y\in B. $ If the coefficients
$\alpha_{\nu}^{\mu}$ are independent of $t $ , $\phi$ satisfies the cocycle
property $\phi_{s}\left(  y\right)  \phi_{r}^{s}\left(  y\right)  =\phi
_{s+r}\left(  y\right)  $ , where $\phi_{t}^{s} $ is the solution to (1) with
$A_{\nu}^{\mu}\left(  t\right)  $ replaced by $A_{\nu}^{s\mu}\left(  t\right)
. $ Define the tensors $a_{\nu}^{\mu}=\alpha_{\nu}^{\mu}\left(  y\right)  $
also for $\mu=+ $ and $\nu=- $ , by
\[
\alpha_{\nu}^{+}\left(  y\right)  =0=\alpha_{-}^{\mu}\left(  y\right)
,\qquad\forall y\in B,
\]
and then one can extend the summation in (\ref{1.1}) to the trace of the
quadratic matrices $\mathbf{a=}\left[  a_{\nu}^{\mu}\right]  $ so it is also
over $\mu=+ $ , and $\nu=- $ . By such an extension the multiplication table
for $\mathrm{d}A\left(  \mathbf{a}\right)  =\mathrm{d}A_{\mu}^{\nu}a_{\nu
}^{\mu}=\boldsymbol{a}\mathrm{d}\boldsymbol{A }$ can be written as
\[
\mathrm{d}A\left(  \mathbf{b}\right)  \mathrm{d}A\left(  \mathbf{a}\right)
=\mathrm{d}A\left(  \mathbf{ba}\right)  ,\quad\mathbf{ba}=\left[  b_{\lambda
}^{\mu}a_{\nu}^{\lambda}\right]
\]
in terms of the usual matrix product $b_{\lambda}^{\mu}a_{\nu}^{\lambda
}=b_{\bullet}^{\mu}a_{\nu}^{\bullet} $ and the involution $\mathbf{a}%
\mapsto\mathbf{a}^{\flat} $ can be obtained by the pseudo-Hermitian
conjugation $a_{\beta}^{\flat\nu}=g^{\nu\kappa}a_{\kappa}^{\mu*}g_{\mu\beta} $
respectively to the indefinite (Minkowski) metric tensor $\mathbf{g}=\left[
g_{\mu\nu}\right]  $ and its inverse $\mathbf{g}^{-1}=\left[  g^{\mu\nu
}\right]  $ , given by $g_{\mu\nu}=\delta_{-\nu}^{\mu}=g^{\mu\nu} $ .

Let us prove that the ''spatial'' part $\boldsymbol{\lambda}=\left(
\lambda_{\nu}^{\mu}\right)  _{\nu\neq-}^{\mu\neq+} $ of the quantum stochastic
germ $\lambda_{\nu}^{\mu}\left(  y\right)  =\delta_{\nu}^{\mu}+\alpha_{\nu
}^{\mu}\left(  y\right)  $ for a PD cocycle exponent $\phi$ must be
conditionally PD in the following sense.

\begin{theorem}
Suppose that the quantum stochastic equation (\ref{1.1}) with $\phi_{0}\left(
y\right)  =y$ has a PD solution in the sense of positive definiteness
(\ref{0.1}) of the matrix $\left[  \phi_{t}\left(  y_{i}^{\star}y_{k}\right)
\right]  $, $\forall t>0$ . Then the germ-matrix $\boldsymbol{\lambda
}=\boldsymbol{p}+\boldsymbol{\alpha}$ to $\boldsymbol{p}=\left(  \delta_{\nu
}^{\mu}\right)  _{\nu\neq-}^{\mu\neq+}$ satisfies the CPD property
\[
\sum_{j}\boldsymbol{e}\boldsymbol{\zeta}_{j}=0\Rightarrow\sum_{i,k}%
\langle\boldsymbol{\zeta}_{i}|\boldsymbol{\lambda}\left(  y_{i}^{\star}%
y_{k}\right)  \boldsymbol{\zeta}_{k}\rangle\geq0.
\]
Here $\boldsymbol{\zeta}\in\mathbb{C\oplus}\mathbb{\ell}_{\mathbb{N}}^{2}$ ,
$\boldsymbol{e}=\left(  e_{\nu}^{\mu}\right)  _{\nu\neq-}^{\mu\neq+}$ ,
$e_{\nu}^{\mu}=\delta_{\nu}^{+}\delta_{-}^{\mu}$ is the one-dimensional
projector, written both with $\boldsymbol{\lambda}$ in the matrix form as
\begin{equation}
\boldsymbol{\lambda}=\left(
\begin{array}
[c]{cc}%
\lambda & \lambda_{\bullet}\\
\lambda^{\bullet} & \lambda_{\bullet}^{\bullet}%
\end{array}
\right)  ,\qquad\boldsymbol{e}=\left(
\begin{array}
[c]{cc}%
1 & 0\\
0 & 0
\end{array}
\right)  ,\qquad\label{1.3}%
\end{equation}
where $\lambda=\alpha_{+}^{-},\quad$ $\lambda^{m}=\alpha_{+}^{m},\quad$
$\lambda_{n}=\alpha_{n}^{-},\quad\lambda_{n}^{m}=\delta_{n}^{m}+\alpha_{n}%
^{m}$ , with $\delta_{n}^{m}\left(  y\right)  =\delta_{n}^{m}$ such that
$\lambda\left(  y^{\star}\right)  =\lambda\left(  y\right)  ^{\dagger},\qquad$
$\lambda^{n}\left(  y^{\star}\right)  =\lambda_{n}\left(  y\right)  ^{\dagger
},\qquad$ $\lambda_{n}^{m}\left(  y^{\star}\right)  =\lambda_{m}^{n}\left(
y\right)  ^{\dagger}.$
\end{theorem}

%

%TCIMACRO{\TeXButton{Proof}{\proof}}%
%BeginExpansion
\proof
%EndExpansion
Let us denote by $\mathcal{D}$ the $\mathbb{C}$ -span $\left\{  \sum_{f}%
\xi^{f}\otimes f^{\otimes}:\xi^{f}\in\mathbb{C},f^{\bullet}\in\mathbb{\ell
}_{\mathbb{N}}^{2}\otimes L^{2}\left(  \mathbb{R}_{+}\right)  \right\}  $ of
coherent (exponential) functions $f^{\otimes}t\left(  \tau\right)
=\bigotimes_{t\in\tau}f^{\bullet}\left(  t\right)  $ , given for each finite
subset $\tau=\left\{  t_{1},...,t_{n}\right\}  \subseteq\mathbb{R}_{+} $ by
tensors $f^{\otimes}\left(  \tau\right)  =f^{n_{1}}\left(  t_{1}\right)
...f^{n_{N}}\left(  t_{N}\right)  $ , where $f^{n},n=\mathbb{N}$ are
square-integrable complex functions on $\mathbb{R}_{+}$ and $\xi^{f}=0$ for
almost all $f^{\bullet}=\left(  f^{n}\right)  $ . The co-isometric shift
$T_{s}$ intertwining $A^{s}\left(  t\right)  $ with $A\left(  t\right)
=T_{s}A^{s}\left(  t\right)  T_{s}^{\dagger}$ is defined on $\mathcal{D}$ by
$T_{s}\left(  f^{\otimes}\right)  \left(  \tau\right)  =f^{\otimes}\left(
\tau+s\right)  $ . The PD property (\ref{0.1}) of the quantum stochastic
adapted map $\phi_{t}$ into the $\mathcal{D}$ -forms $\left\langle \eta
|\phi_{t}\left(  y\right)  \eta\right\rangle $ , for $\eta\in\mathcal{D}$ can
be obviously written as
\begin{equation}
\sum_{i,k}\sum_{f,h}\overline{\xi}_{i}^{f}\phi_{t}\left(  f^{\bullet}%
,y_{i}^{\star}y_{k},h^{\bullet}\right)  \xi_{k}^{h}\geq0,\qquad\label{1.5}%
\end{equation}
for any sequence $y_{j}\in B,j=1,2,...$ , where
\[
\phi_{t}\left(  f^{\bullet},y,h^{\bullet}\right)  =\left\langle f^{\otimes
}|\phi_{t}\left(  y\right)  h^{\otimes}\right\rangle e^{-\int_{t}^{\infty
}f^{\bullet}\left(  s\right)  ^{\dagger}h^{\bullet}\left(  s\right)
\mathrm{d}s},
\]
$\xi^{f}\neq0$ only for a finite subset of $f^{\bullet}\in\left\{
f_{i}^{\bullet},i=1,2,...\right\}  $. If the $\mathcal{D}$ -form $\phi
_{t}\left(  y\right)  $ satisfies the stochastic equation (\ref{1.1}), the
complex function $\phi_{t}\left(  f^{\bullet},y,h^{\bullet}\right)  $
satisfies the differential equation
\begin{align*}
\frac{\mathrm{d}}{\mathrm{d}t}\ln\phi_{t}\left(  f^{\bullet},y,h^{\bullet
}\right)   &  =f^{\bullet}\left(  t\right)  ^{\dagger}h^{\bullet}\left(
t\right)  +\sum_{m,n=1}^{\infty}f^{m}\left(  t\right)  ^{\ast}\alpha_{n}%
^{m}\left(  y\right)  h^{n}\left(  t\right) \\
&  +\sum_{m=1}^{\infty}f^{m}\left(  t\right)  ^{\ast}\alpha_{+}^{m}\left(
y\right)  +\sum_{n=1}^{\infty}\alpha_{n}^{-}\left(  y\right)  h^{n}\left(
t\right)  \phi+\alpha_{+}^{-}\left(  y\right)
\end{align*}
where $f^{\bullet}\left(  t\right)  ^{\dagger}h^{\bullet}\left(  t\right)
=\sum_{n=1}^{\infty}f^{n}\left(  t\right)  ^{\ast}h^{n}\left(  t\right)  $.
The positive definiteness, (\ref{1.5}), ensures the conditional positivity
\[
\sum_{j}\sum_{f}\xi_{j}^{f}=0\Rightarrow\sum_{i,k}\sum_{f,h}\overline{\xi}%
_{i}^{f}\lambda_{t}\left(  f^{\bullet},y_{i}^{\star}y_{k},h^{\bullet}\right)
\xi_{k}^{h}\geq0
\]
of the form $\lambda_{t}\left(  f^{\bullet},y,h^{\bullet}\right)  =\frac{1}%
{t}\left(  \phi_{t}\left(  f^{\bullet},y,h^{\bullet}\right)  -1\right)  $ for
each $t>0$ and any $y_{j}\in B.$ This applies also for the limit $\lambda_{0}$
at $t\downarrow0$ , coinciding with the quadratic form
\[
\frac{\mathrm{d}}{\mathrm{d}t}\phi_{t}\left(  f^{\bullet},y,h^{\bullet
}\right)  |_{t=0}=\sum_{m,n}\bar{a}^{m}\lambda_{n}^{m}\left(  y\right)
c^{n}+\sum_{m}\bar{a}^{m}\lambda^{m}\left(  y\right)  +\sum_{n}\lambda
_{n}\left(  y\right)  c^{n}+\lambda\left(  y\right)  ,
\]
where $a^{\bullet}=f^{\bullet}\left(  0\right)  ,\quad c^{\bullet}=h^{\bullet
}\left(  0\right)  $ , and the $\lambda$ 's are defined in (\ref{1.3}). Hence
the form
\[
\sum_{i,k}\sum_{\mu,\nu}\overline{\zeta}_{i}^{\mu}\lambda_{\nu}^{\mu}\left(
y_{i}^{\star}y_{k}\right)  \zeta_{k}^{\nu}:=\sum_{i,k}\overline{\zeta}%
_{i}\lambda\left(  y_{i}^{\star}y_{k}\right)  \zeta_{k}%
\]%
\[
+\sum_{i,k}\left(  \sum_{n}\overline{\zeta}_{i}\lambda_{n}\left(  y_{i}%
^{\star}y_{k}\right)  \zeta_{k}^{n}+\sum_{m}\overline{\zeta}_{i}^{m}%
\lambda^{m}\left(  y_{i}^{\star}y_{k}\right)  \zeta_{k}+\sum_{m,n}%
\overline{\zeta}_{i}^{m}\lambda_{n}^{m}\left(  y_{i}^{\star}y_{k}\right)
\zeta_{k}^{n}\right)
\]
with $\zeta=\sum_{f}\xi^{f},\quad\zeta^{\bullet}=\sum_{f}\xi^{f}a_{f}%
^{\bullet}$ , where $a_{f}^{\bullet}=f^{\bullet}\left(  0\right)  $ , is
positive if $\sum_{j}\zeta_{j}=0.$ The components $\zeta$ and $\zeta^{\bullet
}$ of these vectors are independent because for any $\zeta\in\mathbb{C}$ and
$\zeta^{\bullet}=\left(  \zeta^{1},\zeta^{2},...\right)  \in\mathbb{\ell
}_{\mathbb{N}}^{2}$ there exists such a function $a^{\bullet}\mapsto\xi^{a}$
on $\mathbb{\ell}_{\mathbb{N}}^{2}$ with a finite support, that $\sum_{a}%
\xi^{a}=\zeta,\quad\sum_{a}\xi^{a}a^{\bullet}=\zeta^{\bullet}$ , namely,
$\xi^{a}=0$ for all $a^{\bullet}\in\mathbb{\ell}_{\mathbb{N}}^{2}$ except
$a^{\bullet}=0$ , for which $\xi^{a}=\zeta-\sum_{n=1}^{\infty}\zeta^{n}$ and
$a^{\bullet}=e_{n}^{\bullet}$ , the $n$ -th basis element in $\mathbb{\ell
}_{\mathbb{N}}^{2}$ , for which $\xi^{a}=\zeta^{n}.$ This proves the complete
positivity of the matrix form $\boldsymbol{\lambda}$ , with respect to the
matrix orthoprojector $\boldsymbol{p}_{0}$ defined in (\ref{1.3}) on the
ket-vectors $\boldsymbol{\zeta}=\left(  \zeta^{\mu}\right)  $
%TCIMACRO{\TeXButton{End Proof}{\endproof}}%
%BeginExpansion
\endproof
%EndExpansion

\section{A Dilation Theorem for the Form-Generator.}

The CPD property of the germ-matrix $\boldsymbol{\lambda}$ with respect to the
projective matrix $\boldsymbol{p}_{0}$ (\ref{1.3}) obviously implies the
positivity of the dissipation form
\begin{equation}
\sum_{x,z}\left\langle \boldsymbol{\zeta}^{x}|\boldsymbol{\Delta}\left(
x,z\right)  \boldsymbol{\zeta}^{z}\right\rangle :=\sum_{k,l}\sum_{\mu,\nu
}\left\langle \zeta_{k}^{\mu}|\Delta_{\nu}^{\mu}\left(  y_{k},y_{l}\right)
\zeta_{l}^{\nu}\right\rangle , \label{2.1}%
\end{equation}
$\qquad$ where $\zeta^{-}=\zeta=\zeta^{+}$ and $\zeta_{j}=\zeta^{y_{j}} $ for
any (finite) sequence $y_{j}\in B$ , $j=1,2,...$ , corresponding to non-zero
$\boldsymbol{\zeta}_{y}\in\mathbb{C\oplus\ell}_{\mathbb{N}}^{2}$ . Here
$\boldsymbol{\Delta}=\left(  \Delta_{\nu}^{\mu}\right)  _{\nu=+,\bullet}%
^{\mu=-,\bullet}$ is the stochastic dissipator
\[
\boldsymbol{\Delta}\left(  x,z\right)  =\boldsymbol{\lambda}\left(  x^{\star
}z\right)  -\boldsymbol{e}\boldsymbol{\lambda}\left(  z\right)
-\boldsymbol{\lambda}\left(  x^{\star}\right)  \boldsymbol{e}+\boldsymbol{e}%
\boldsymbol{\lambda}\left(  1\right)  \boldsymbol{e }%
\]
with the elements
\begin{align}
\Delta_{n}^{m}\left(  x,z\right)   &  =\alpha_{n}^{m}\left(  x^{\star
}z\right)  +\delta_{n}^{m},\label{2.2}\\
\Delta_{n}^{-}\left(  x,z\right)   &  =\alpha_{n}^{-}\left(  x^{\star
}z\right)  -\alpha_{n}^{-}\left(  z\right)  =\Delta_{+}^{n}\left(  z,x\right)
^{\dagger},\nonumber\\
\Delta_{+}^{-}\left(  x,z\right)   &  =\alpha_{+}^{-}\left(  x^{\star
}z\right)  -\alpha_{+}^{-}\left(  z\right)  -\alpha_{+}^{-}\left(  x^{\star
}\right)  +d,\nonumber
\end{align}
where $d=\alpha_{+}^{-}\left(  1\right)  \leq0$ ( $d=0$ for the case of the
martingale $M_{t}=\phi_{t}\left(  1\right)  $ ). In particular the
matrix-valued map $\lambda_{\bullet}^{\bullet}=\left[  \lambda_{n}^{m}\right]
$ is PD. If the functions $\lambda^{m}$ , $\lambda_{n},\lambda$ have the form
\begin{equation}
\lambda^{m}\left(  y\right)  =\varphi^{m}\left(  y\right)  -c^{m},\quad
\lambda_{n}\left(  y\right)  =\varphi_{n}\left(  y\right)  -c_{n},\quad
\lambda\left(  y\right)  =\varphi\left(  y\right)  -c \label{2.3}%
\end{equation}
such that $\boldsymbol{\varphi}=\boldsymbol{\lambda}-\boldsymbol{c}$ , is a PD
map for a constant Hermitian matrix $\boldsymbol{c}=\left(  c_{\nu}^{\mu
}\right)  _{\nu\neq-}^{\mu\neq+}$ , the CPD condition is fulfilled for
$\boldsymbol{\lambda}$ .

In order to make the formulation of the following dilation theorem as concise
as possible, we need the notion of the $\flat$ -representation of $B $ in a
pseudo-Hilbert space $\mathcal{E}=\mathbb{C}\oplus\mathcal{K}\oplus\mathbb{C}$
with respect to the indefinite metric
\begin{equation}
\left(  \xi|\xi\right)  =2\operatorname{Re}\overline{\xi}^{-}\xi^{+}+\left\|
\xi^{\circ}\right\|  ^{2}+|\xi^{+}|^{2}d \label{2.4}%
\end{equation}
for the triples $\xi=\left(  \xi^{-},\xi^{\circ},\xi^{+}\right)
\in\mathcal{E}$ , where $\xi^{-},\xi^{+}\in\mathbb{C},\quad\xi^{\circ}%
\in\mathcal{K},\quad\mathcal{K}$ is a pre-Hilbert space. The operators $A$ in
this space are given by the $3\times3$ -block-matrices $\mathbf{A}=\left[
A_{\nu}^{\mu}\right]  _{\nu=+,\circ,+}^{\mu=-,\circ,+}$ , and the
pseudo-Hermitian conjugation $\left(  A^{\flat}\xi|\xi\right)  =\left(
\xi|A\xi\right)  $ is given by the usual Hermitian conjugation $A_{\nu
}^{\dagger\mu}=A_{\mu}^{\nu\ast}$ as $\mathbf{A}^{\flat}=\mathbf{G}%
^{-1}\mathbf{A}^{\dagger}\mathbf{G}$ respectively to the indefinite metric
tensor $\mathbf{G}=\left[  G_{\mu\nu}\right]  $ and its inverse $\mathbf{G}%
^{-1}=\left[  G^{\mu\nu}\right]  $ , given by
\begin{equation}
\mathbf{G}=\left[
\begin{array}
[c]{ccc}%
0 & 0 & 1\\
0 & I_{\circ}^{\circ} & 0\\
1 & 0 & d
\end{array}
\right]  ,\qquad\mathbf{G}^{-1}=\left[
\begin{array}
[c]{ccc}%
-d & 0 & 1\\
0 & I_{\circ}^{\circ} & 0\\
1 & 0 & 0
\end{array}
\right]  \label{2.5}%
\end{equation}
with a real $d$ , where $I_{\circ}^{\circ}$ is the identity operator in
$\mathcal{K}$ . The algebras of all operators $A$ on $\mathcal{K}$ and
$\mathcal{E}$ with $A^{\dagger}\mathcal{K}\subseteq\mathcal{K}$ and $A^{\flat
}\mathcal{E}\subseteq\mathcal{E}$ are denoted by $\mathcal{A}\left(
\mathcal{K}\right)  $ and $\mathcal{A}\left(  \mathcal{E}\right)  $ .

\begin{theorem}
The following are equivalent:

\begin{enumerate}
\item The dissipator (\ref{2.2}), defined by the $\flat$ -map $\alpha$ with
$\alpha_{+}^{-}\left(  1\right)  =d$ , is positive definite:
\[
\sum_{x,z}\left\langle \boldsymbol{\zeta}_{x}|\boldsymbol{\Delta}\left(
x,z\right)  \boldsymbol{\zeta}_{z}\right\rangle \geq0
\]

\item There exist: a pre-Hilbert space $\mathcal{K}$ , a unital $\dagger$ -
representation $j$ in $\mathcal{A}\left(  \mathcal{K}\right)  $ ,
\begin{equation}
j\left(  x^{\star}z\right)  =j\left(  x\right)  ^{\dagger}j\left(  z\right)
,\quad j\left(  1\right)  =I, \label{2.6}%
\end{equation}
of the $\star$ -multiplication structure of $B$ , a $j$ -cocycle on $B$ ,
\begin{equation}
k\left(  x^{\star}z\right)  =j\left(  x\right)  ^{\dagger}k\left(  z\right)
+k\left(  x^{\star}\right)  , \label{2.7}%
\end{equation}
having values in $\mathcal{K}$ , and a function $l:B\rightarrow\mathbb{C}$ ,
having the coboundary property
\begin{equation}
l\left(  x^{\star}z\right)  =l\left(  z\right)  +l\left(  x^{\star}\right)
+k^{*}\left(  x^{\star}\right)  k\left(  z\right)  , \label{2.8}%
\end{equation}
with $k^{*}\left(  y^{\star}\right)  =k\left(  y\right)  ^{*},l\left(
y^{\star}\right)  =l\left(  y\right)  ^{*}$ $,\quad$ such that $\lambda\left(
y\right)  =l\left(  y\right)  +d$ ,
\[
\lambda_{n}\left(  y^{\star}\right)  =k\left(  y\right)  ^{\dagger}%
L_{n}^{\circ}+L_{n}^{-}=\lambda^{n}\left(  y\right)  ^{\dagger},
\]
and $\lambda_{n}^{m}\left(  y\right)  =L_{m}^{\circ*}j\left(  y\right)
L_{n}^{\circ}$ for some elements $L_{n}^{\circ}\in\mathcal{K}$ with the
adjoints $L_{n}^{\circ*}=L_{\circ}^{n}:\mathcal{K}\rightarrow\mathbb{C}$ and
$L_{n}^{-}\in\mathbb{C}$ .

\item There exist a pseudo-Hilbert space, $\mathcal{E}$ , namely,
$\mathbb{C}\oplus\mathcal{K}\oplus\mathbb{C}$ with the indefinite metric
tensor $\mathbf{G}=\left[  G_{\mu\nu}\right]  $ given above for $\mu
,\nu=-,\circ,+ $ , and $d=\lambda\left(  1\right)  $ , a unital $\flat$
-representation $\mathbf{\jmath}=\left[  \jmath_{\nu}^{\mu}\right]
_{\nu=-,\circ,+}^{\mu=-,\circ,+} $ of the $\star$ -multiplication structure of
$B$ on $\mathcal{E}$ :
\begin{equation}
\mathbf{\jmath}\left(  x^{\star}z\right)  =\mathbf{\jmath}\left(  x\right)
^{\flat}\mathbf{\jmath}\left(  z\right)  ,\quad\mathbf{\jmath}\left(
1\right)  =\mathbf{I} \label{2.9}%
\end{equation}
with $\mathbf{\jmath}\left(  y\right)  ^{\flat}=\mathbf{G}^{-1}\mathbf{\jmath
}\left(  y\right)  ^{\dagger}\mathbf{G}$ , given by the matrix elements
\[
\jmath_{\circ}^{\circ}=j,\quad\jmath_{+}^{\circ}=k,\quad\jmath_{\circ}%
^{-}=k^{*},\quad\jmath_{+}^{-}=l,\quad\jmath_{-}^{-}=1=\jmath_{+}^{+}%
\]
and all other $\jmath_{\nu}^{\mu}=0$ , and a linear operator $\mathbf{L}%
:\mathbb{C}\oplus\ell_{\mathbb{N}}^{2}\rightarrow\mathcal{E}$ , with the
components $\left[  L^{\mu},L_{\bullet}^{\mu}\right]  $ , where
\[
L^{-}=0,\quad L^{\circ}=0,\quad L^{+}=1,\quad L_{\bullet}^{-}=\left(
L_{n}^{-}\right)  ,\quad L_{\bullet}^{\circ}=\left(  L_{n}^{\circ}\right)
,\quad L_{\bullet}^{+}=0,
\]
and $\mathbf{L}^{\flat}=\left(
\begin{array}
[c]{ccc}%
1 & 0 & \delta\\
0 & L_{\circ}^{\bullet} & L_{+}^{\bullet}%
\end{array}
\right)  =\mathbf{L}^{\dagger}\mathbf{G}$ , where $L_{\circ}^{\bullet
}=L_{\bullet}^{\circ\dagger},L_{+}^{\bullet}=L_{\bullet}^{-\dagger}$ , such
that
\begin{equation}
\mathbf{L}^{\flat}\mathbf{\jmath}\left(  y\right)  \mathbf{L}%
=\boldsymbol{\lambda}\left(  y\right)  ,\qquad\forall y\in B. \label{2.10}%
\end{equation}

\item The germ-matrix $\boldsymbol{\lambda}\left(  y\right)  =\left(
\alpha_{\nu}^{\mu}\left(  y\right)  +\delta_{\nu}^{\mu}\right)  _{\nu\neq
-}^{\mu\neq+}$ is CPD with respect to the orthoprojector $\boldsymbol{e}$ ,
defined in ( $\ref{1.3})$ :
\[
\sum_{y}\boldsymbol{e}\boldsymbol{\zeta}^{y}=0\Rightarrow\sum_{x,z}%
\langle\boldsymbol{\zeta}^{x}|\boldsymbol{\lambda}\left(  x^{\star}z\right)
\boldsymbol{\zeta}^{z}\rangle\geq0.
\]

\end{enumerate}
\end{theorem}

%

%TCIMACRO{\TeXButton{Proof}{\proof}}%
%BeginExpansion
\proof
%EndExpansion
Similar to the dilation theorem in \cite{4}, see also \cite{6}, \cite{7},
\cite{8}%
%TCIMACRO{\TeXButton{End Proof}{\endproof}}%
%BeginExpansion
\endproof
%EndExpansion

\section{Pseudo-Poisson processes and their generators.}

Let us consider the case $B=1\oplus\mathfrak{b} $ of the unital semigroup for
a $\star$ -algebra $\mathfrak{b} $ with $\boldsymbol{\lambda}\left(  1\oplus
b\right)  =\boldsymbol{d}+\boldsymbol{\gamma}\left(  b\right)  $ given by a
linear matrix -function
\[
\boldsymbol{\gamma}=\left(
\begin{array}
[c]{cc}%
\gamma & \gamma_{\bullet}\\
\gamma^{\bullet} & \gamma_{\bullet}^{\bullet}%
\end{array}
\right)  =\boldsymbol{\lambda}-\boldsymbol{d},\quad\boldsymbol{d}=\left(
\begin{array}
[c]{cc}%
d & d_{\bullet}\\
d^{\bullet} & d_{\bullet}^{\bullet}%
\end{array}
\right)  =\boldsymbol{\lambda}\left(  1\right)
\]
of $b\in\mathfrak{b} $ for $y=1\oplus b $ . Following \cite{4}, the linear
quantum stochastic process $\Lambda\left(  t\right)  :b\mapsto
\boldsymbol{\gamma}\left(  b\right)  \boldsymbol{A}\left(  t\right)  $ with
independent increments, generating together with $A\left(  t,\mathbf{d}%
\right)  =A_{\mu}^{\nu}\left(  t\right)  d_{\nu}^{\mu}$ the stochastic PD
exponent
\[
\phi_{t}\left(  1\oplus b\right)  =:\exp\left[  A\left(  t,\mathbf{d}\right)
+\Lambda\left(  t,b\right)  \right]  :\quad b\in\mathfrak{b}%
\]
as the solution of the equation (\ref{1.1}), will be called the
pseudo-Poissonian\cite{4} over the algebra $\mathfrak{b} $ .

If $B$ is a unit ball of an operator algebra $\mathcal{B}$ , the linear
form-generator can be extended to the whole algebra. The structure (\ref{2.3})
of the linear form-generator for PD cocycles over an operator algebra
$\mathcal{B}$ is a consequence of the cocycle equation (\ref{2.7}), according
to which $j\left(  0\right)  k\left(  y\right)  =0$ , where
\begin{equation}
k\left(  y\right)  =j\left(  y\right)  \varsigma-\varsigma,.\quad
\varsigma=-k\left(  0\right)  . \label{3.1}%
\end{equation}
Denoting by $\varsigma^{\dagger}$ the linear functional $\xi^{\circ}%
\mapsto\left(  \varsigma|\xi^{\circ}\right)  $ on $\mathcal{K}$ corresponding
to the $\varsigma\in\mathcal{K}$ , the condition (\ref{2.8}) yields
\begin{equation}
l\left(  y\right)  =\frac{1}{2}\left(  \varsigma^{\dagger}k\left(  y\right)
+k^{\ast}\left(  y\right)  \varsigma\right)  =\varsigma^{\dagger}j\left(
y\right)  \varsigma-\varsigma^{\dagger}\varsigma. \label{3.2}%
\end{equation}
Hence, in addition to $\lambda_{n}^{m}\left(  y\right)  =L_{m}^{\circ\dagger
}j\left(  y\right)  L_{n}^{\circ}$ one can obtain the structure (\ref{2.3})
with
\begin{equation}
\varphi\left(  y\right)  =\varsigma^{\dagger}j\left(  y\right)  \varsigma
,\quad\varphi_{n}\left(  y\right)  =\varsigma^{\dagger}j\left(  y\right)
L_{n}^{\circ},\quad\varphi^{m}\left(  y\right)  =L_{m}^{\circ\dagger}j\left(
y\right)  \varsigma,\quad\label{3.3}%
\end{equation}
and $\kappa=\varsigma^{\dagger}\varsigma-\delta$ , $\kappa_{n}=\varsigma
^{\dagger}L_{n}^{\circ}-L_{n}^{-}$ . Thus, $\boldsymbol{\lambda}\left(
y\right)  =\boldsymbol{\varphi}\left(  y\right)  -\boldsymbol{\kappa}$ , where
$\boldsymbol{\varphi}$ is a completely positive nonlinear map of $B$ into the
space $\mathcal{M}\left(  \mathbb{C\oplus\ell}_{\mathbb{N}}^{2}\right)  $ of
complex matrices $\boldsymbol{x}=\left(  x_{\nu}^{\mu}\right)  $ . Moreover,
$\boldsymbol{\varphi}$ is uniquely defined as the birth-map by the condition
$\boldsymbol{\varphi}\left(  0\right)  =0$ with $\boldsymbol{\kappa
}=-\boldsymbol{\lambda}\left(  0\right)  =\left(  \kappa_{\nu}^{\mu}\right)
$, where $\kappa_{+}^{-}=\kappa,\kappa_{n}^{-}=\kappa_{n},\kappa_{+}%
^{m}=\overline{\kappa}_{m}$ , and $\kappa_{n}^{m}=-\lambda_{\nu}^{\mu}\left(
0\right)  $ , constituting a negative-definite matrix $\kappa_{\bullet
}^{\bullet}=\left[  \kappa_{n}^{m}\right]  $ . Any germ-matrix
$\boldsymbol{\lambda}$ whose components are decomposed into the sums of the
components $\varphi_{\nu}^{\mu}\quad$ of a PD map $\boldsymbol{\varphi}$ and
$\boldsymbol{\lambda}\left(  0\right)  $ , are obviously CPD with respect to
the orthoprojector $\boldsymbol{p}_{0}$ in (\ref{1.5}). As follows from the
dilation theorem, there exists a family $\varsigma_{-}=\varsigma=\varsigma
_{+},\quad\varsigma_{n}=L_{n}^{\circ}-j\left(  0\right)  L_{n}^{\circ},\quad
n\in\mathbb{N}$ of vectors $\varsigma_{\nu}\in\mathcal{K}$ with $j\left(
0\right)  \varsigma_{\nu}=0$ such that $\varphi_{\nu}^{\mu}\left(  y\right)
=\varsigma_{\mu}^{\dagger}j\left(  y\right)  \varsigma_{\nu}$ for all $\mu
\in\left\{  -,1,2,...\right\}  $ , $\nu\in\left\{  +,1,2,...\right\}  $. Thus
the equation (\ref{1.1}) for a completely positive exponential cocycle with
bounded stochastic derivatives has the following general form
\[
\mathrm{d}\phi_{t}\left(  y\right)  +\left(  \gamma-\varsigma^{\dagger
}j\left(  y\right)  \varsigma\right)  \phi_{t}\left(  y\right)  \mathrm{d}%
t=\sum_{m,n=1}^{\infty}\left(  \varsigma_{m}^{\dagger}j\left(  y\right)
\varsigma_{n}-\gamma_{n}^{m}\right)  \phi_{t}\left(  y\right)  \mathrm{d}%
A_{m}^{n}%
\]%
\begin{equation}
+\sum_{m=1}^{\infty}\left(  \varsigma_{m}^{\dagger}j\left(  y\right)
\varsigma-\gamma_{m}^{\dagger}\right)  \phi_{t}\left(  y\right)
\mathrm{d}A_{m}^{+}+\sum_{n=1}^{\infty}\left(  \varsigma^{\dagger}j\left(
y\right)  \varsigma_{n}-\gamma_{n}\right)  \phi_{t}\left(  y\right)
\mathrm{d}A_{-}^{n},\qquad\label{3.4}%
\end{equation}
where $\gamma_{\nu}^{\mu}=-\alpha_{\nu}^{\mu}\left(  0\right)  $ . If
$M_{t}=\phi_{t}\left(  1\right)  $ is a martingale, the normalization
condition $\sum_{k=1}^{\infty}\varsigma^{k\dagger}\varsigma^{k}=\kappa$ (
$\leq\kappa$ if submartingale).\allowbreak

In the particular case $\mathcal{K}=\mathbb{C\oplus}\mathfrak{h} $ , $j\left(
y\right)  =1\oplus y $ , where $\mathfrak{h} $ is a Hilbert space of a
representation $\mathcal{B}\subseteq\mathcal{B}\left(  \mathfrak{h}\right)  $
of the C*-algebra $\mathcal{B} $ in the operator algebra $\mathcal{B}\left(
\mathfrak{h}\right)  $ , this gives a quantum stochastic generalization of the
Poissonian birth semigroups \cite{1} with the affine generators $\alpha_{\nu
}^{\mu}\left(  y\right)  =\varsigma_{\mu}^{\dagger}X\varsigma_{\nu}%
-\gamma_{\nu}^{\mu}$ . In the more general case when the space $\mathcal{K} $
is embedded into the Hilbert sum of all tensor powers of the space
$\mathfrak{h} $ such that $j\left(  y\right)  =\oplus_{k=0}^{\infty}y^{\otimes
k} $ , the birth function $\boldsymbol{\varphi}$ is described by the
components
\begin{align}
\varphi_{n}^{m}\left(  y\right)   &  =\sum_{k=0}^{\infty}\varsigma
_{m}^{k\dagger}y^{\otimes k}\varsigma_{n}^{k},\qquad\varphi\left(  y\right)
=\sum_{k=1}^{\infty}\varsigma^{k\dagger}y^{\otimes k}\varsigma^{k}%
\label{3.5}\\
\varphi^{m}\left(  y\right)   &  =\sum_{k=1}^{\infty}\varsigma_{m}^{k\dagger
}y^{\otimes k}\varsigma^{k},\qquad\varphi_{n}\left(  y\right)  =\sum
_{k=1}^{\infty}\varsigma^{k\dagger}y^{\otimes k}\varsigma_{n}^{k}%
\qquad\nonumber
\end{align}
with $\varsigma^{k},\varsigma_{n}^{k}\in\mathfrak{h}^{\otimes k} $ .

Note, if $\mathcal{B} $ is a W*-algebra and the germ map $\boldsymbol{\lambda
}$ is w*-analytic, the completely positive function $\boldsymbol{\varphi}$ is
also analytic, being defined by a w*-analytical representation $j=\oplus
_{k=0}^{\infty}i^{\otimes k} $ in a full Fock space $\mathcal{K}=\oplus
_{k=0}^{\infty}\mathcal{H}^{\otimes k} $ , where $i $ is a (linear)
w*-representation of $\mathcal{B} $ on a Hilbert space $\mathcal{H} $ . This
gives the general form for the w*-analytical quantum stochastic quasi-Poisson
birth process over the algebra $\mathcal{B} $ .

The next theorem proves that these structural conditions which are necessary
for complete positivity of the stochastic exponents, given by the equation
(\ref{1.1}), are also sufficient. In particular it proves the existence of the
quantum birth cocycle $\phi$ for a given generating stochastic birth
matrix-function $\boldsymbol{\varphi}. $

\begin{theorem}
Let the structural maps $\boldsymbol{\lambda}$ of the quantum stochastic PD
exponent $\phi$ over the unit ball of an operator algebra $\mathcal{B}$ . Then
they are bounded in the unit ball of $\mathcal{B}$ ,
\[
||\lambda||<\infty,\qquad||\lambda_{\bullet}||=\left(  \sum_{n=1}^{\infty
}\left\|  \lambda_{n}\right\|  ^{2}\right)  ^{\frac12}=\left\|  \lambda
^{\bullet}\right\|  <\infty,\qquad\left\|  \lambda_{\bullet}^{\bullet
}\right\|  =\left\|  \lambda_{\bullet}^{\bullet}\left(  1\right)  \right\|
<\infty,
\]
where $\left\|  \lambda\right\|  =\sup\left\{  \left\|  \lambda\left(
y\right)  \right\|  :\left\|  y\right\|  <1\right\}  ,\left\|  \lambda
_{\bullet}^{\bullet}\left(  1\right)  \right\|  =\sup\left\{  \left\langle
\zeta^{\bullet}|\lambda_{\bullet}^{\bullet}\left(  1\right)  \zeta^{\bullet
}\right\rangle \left|  \left\|  \zeta^{\bullet}\right\|  <1\right.  \right\}
$ , and have the form (\ref{3.3}) written as
\[
\boldsymbol{\lambda}\left(  y\right)  =\boldsymbol{\varphi}\left(  y\right)
-\boldsymbol{\kappa}%
\]
with $\varphi=\varphi_{+}^{-},\quad\varphi^{m}=\varphi_{+}^{m},\quad
\varphi_{n}=\varphi_{n}^{-}$ and $\varphi_{n}^{m}=\lambda_{n}^{m}$ , composing
a bounded PD map
\begin{equation}
\boldsymbol{\varphi}=\left[
\begin{array}
[c]{cc}%
\varphi & \varphi_{\bullet}\\
\varphi^{\bullet} & \varphi_{\bullet}^{\bullet}%
\end{array}
\right]  ,\quad and\quad\boldsymbol{\kappa}=\left[
\begin{array}
[c]{cc}%
\kappa & \kappa_{\bullet}\\
\kappa_{\bullet}^{*} & 0
\end{array}
\right]  \label{3.6}%
\end{equation}
with arbitrary $\kappa$ and $\kappa_{\bullet}=\left(  \kappa_{1}%
,\kappa_{2,...}\right)  $ . The equation (\ref{3.4}) has the unique PD
solution
\begin{equation}
\phi_{t}\left(  y\right)  =V_{t}^{\dagger}\exp\left[  A_{\bullet}^{+}\left(
t\right)  \varphi^{\bullet}\left(  y\right)  \right]  \varphi_{\bullet
}^{\bullet}\left(  y\right)  ^{A_{\bullet}^{\bullet}\left(  t\right)  }%
\exp\left[  \varphi_{\bullet}\left(  y\right)  A_{-}^{\bullet}\left(
t\right)  \right]  V_{t}\exp\left[  t\varphi\left(  y\right)  \right]  ,
\label{3.7}%
\end{equation}
where $V_{t}=\exp\left[  -\kappa_{\bullet}A_{-}^{\bullet}\left(  t\right)
-\frac12\kappa tI\right]  .$
\end{theorem}

%

%TCIMACRO{\TeXButton{Proof}{\proof}}%
%BeginExpansion
\proof
%EndExpansion
(Sketch) The PD solution to the quantum stochastic equation (\ref{3.4}) can be
obtained by the iteration of the equivalent quantum stochastic integral
equation
\[
\phi_{t}\left(  y\right)  =V_{t}^{\dagger}V_{t}+\int_{0}^{t}V_{s}^{\dagger
}\phi_{t-s}^{s}\left(  y\right)  V_{s}\beta_{\nu}^{\mu}\left(  y\right)
\mathrm{d}A_{\mu}^{\nu}\left(  s\right)
\]
where $\beta_{\nu}^{\mu}\left(  y\right)  =\varphi_{\nu}^{\mu}\left(
y\right)  -\delta_{\nu}^{\mu}$ .Here $V_{t}$ is the exponential vector cocycle
$V_{r}^{s}V_{s}=V_{r+s}$ , resolving the quantum stochastic differential
equation
\[
\mathrm{d}V_{t}+\kappa V_{t}\mathrm{d}t+\sum_{n=1}^{\infty}\kappa_{n}%
V_{t}\mathrm{d}A_{-}^{n}=0
\]
with the initial condition $V_{0}=I$ in $\mathcal{D}$ and with $V_{r}%
^{s}=T_{r}^{\dagger}V_{r}T_{s}$ , shifted by the time-shift co-isometry
$T_{s}$ in $\mathcal{D}$.

\end{document}